\documentclass[11pt]{article}

\usepackage{latexsym}
\usepackage{amssymb}

\newtheorem{theorem}{Theorem}

\newtheorem{lemma}{Lemma}
\newtheorem{remark}{Remark}
\newtheorem{example}{Example}

\begin{document}
{
\begin{center}
{\Large\bf
On the similarity of complex symmetric operators to perturbations of restrictions of normal operators.}
\end{center}
\begin{center}
{\bf Sergey M. Zagorodnyuk}
\end{center}

\noindent
\textbf{Abstract.}
In this paper we consider a problem of the similarity of complex symmetric operators to perturbations of restrictions of
normal operators.
For a subclass of cyclic complex symmetric operators in a finite-dimensional Hilbert space we prove the similarity to
rank-one perturbations of restrictions of normal operators. The main tools are a truncated moment problem in $\mathbb{C}$, and
some objects similar to objects from the theory of spectral problems for Jacobi matrices.

\section{Introduction.}

During past 15 years an increasing interest was devoted to complex symmetric operators and other types of
operators related to a conjugation in a separable Hilbert space $H$, see~\cite{cit5000}, \cite{cit_Guo_Ji_Zhu_2015}, \cite{cit_9000_Zhu_2016},
\cite{Wang_Zhu_2020_Reducing} and papers cited therein.
The conjugation $J$ is an {\it antilinear} operator in $H$ such that $J^2 x = x$, $x\in H$,
and
$$ (Jx,Jy)_H = (y,x)_H,\qquad x,y\in H. $$
Denote
\begin{equation}
\label{f1_1}
[x,y]_J := (x,Jy)_H,\qquad x,y\in H,
\end{equation}
where $(\cdot,\cdot)_H$ is the inner product in a Hilbert space $H$.
Recall that a linear operator $A$ in $H$ is said to be \textit{$J$-symmetric} if
\begin{equation}
\label{f1_2}
[Ax,y]_J = [x,Ay]_J,\qquad x,y\in D(A).
\end{equation}
Observe that for a \textit{bounded} linear operator $A$ condition~(\ref{f1_2}) is
equivalent to the following condition:
\begin{equation}
\label{f1_4}
JAJ = A^*.
\end{equation}
If a linear bounded operator $A$ on a whole Hilbert space $H$ is $J$-symmetric 
for some conjugation $J$ in $H$, then $A$ is said to be \textit{complex symmetric}. The latter notion was introduced
by Garcia and Putinar in~2006 in~\cite{cit_1000_GP}.
It should be noticed that $J$-symmetric operators appeared much earlier, namely, in a paper of Glazman in~1957~\cite{cit100}.
A brief survey of the history of such operators can be found in the introduction of a paper~\cite{cit_3000_Z}.
Garcia and Putinar in~\cite{cit_1000_GP} did not explain their reasons for their notion. Probably, they intended to
have a class of operators which preserve unitary equivalent operators inside the class.

Let $d$ be a fixed integer greater than $1$. We shall say that a matrix $\mathcal{M} = (m_{k,l})_{k,l=0}^{d-1}$, $m_{k,l}\in \mathbb{C}$,  
belongs to the class $\mathfrak{M}_{d;3}^+$, if and only if
the following conditions hold:
\begin{equation}
\label{f1_10}
m_{k,l} = 0,\qquad  k,l\in \mathbb{Z}_{0,d-1}:\ |k-l|>1;
\end{equation}
\begin{equation}
\label{f1_20}
m_{k,l} = m_{l,k},\qquad k,l\in \mathbb{Z}_{0,d-1};
\end{equation}
\begin{equation}
\label{f1_30}
m_{k,k+1} \not= 0,\qquad k\in \mathbb{Z}_{0,d-2}.
\end{equation}

Let $A$ be a linear operator in a \textit{finite-dimensional} Hilbert space $H$ of dimension $d$.
We shall say that $A$ belongs to the class $C_+ = C_+(H)$ if and only if there exists an
orthonormal basis $\{ u_k \}_{k=0}^{d-1}$ in $H$ such that the matrix
\begin{equation}
\label{f1_35}
\mathcal{M} = ( (A u_l,u_k) )_{k,l=0}^{d-1},
\end{equation}
belongs to $\mathfrak{M}_{d;3}^+$.
The above notions of classes $\mathfrak{M}_{d;3}^+$ and $C_+ = C_+(H)$ are similar to notions for complex
symmetric operators in an infinite-dimensional Hilbert space from~\cite{cit_3500_Z}.
Thus, $C_+(H)$ consists of operators in $H$ having a three-diagonal complex symmetric matrix, with non-zero entries on the
first sub-diagonal, with respect to an orthonormal basis of $H$.
We shall characterize these operators below (see Theorem~\ref{t2_1}). 
They form a subclass of cyclic complex symmetric operators in $H$.
Recall that a bounded linear operator $A$ in a Hilbert space $H$ is said to be cyclic if
there exists a vector $x_0\in H$ (cyclic vector) such that
$$ \overline{ \mathop{\rm Lin}\nolimits \{ A^k x_0,\ k\in\mathbb{Z}_+ \} } = H. $$

The main objective of this paper is to show that operators from $C_+(H)$ are similar to rank-one perturbations of restrictions of normal
operators (notice that such restrictions need not to be subnormal operators). For this purpose we shall use the following moment problem:
find a (non-negative) measure $\mu$ on $\mathfrak{B}(\mathbb{C})$
such that
\begin{equation}
\label{f1_50}
\int_{\mathbb{C}} z^k d\mu(z) = s_k,\qquad  k\in\mathbb{Z}_{0,\rho}.
\end{equation}
Here $\{ s_k \}_{k\in\mathbb{Z}_{0,\rho}}$ is a prescribed set of complex numbers (moments); $\rho\in\mathbb{N}$.
A more general problem was recently considered in~\cite{cit_3700_Z} (here we shall not need any results from~\cite{cit_3700_Z}). 
We shall present a transparent construction of atomic solutions for a solvable moment problem~(\ref{f1_50}),
which have arbitrarily big number of atoms.
We shall also need some objects and results similar to those objects and results from the theory of spectral problems for Jacobi matrices.
It is interesting that for these results the complex symmetry was crucial.
Finally, we state some open problems which appear in a natural way from our discussion.

\noindent
{\bf Notations. }
Throughout the whole paper $d$ means a fixed integer greater than $1$.
As usual, we denote by $\mathbb{R}, \mathbb{C}, \mathbb{N}, \mathbb{Z}, \mathbb{Z}_+$
the sets of real numbers, complex numbers, positive integers, integers and non-negative integers,
respectively; $\mathbb{T} = \{ z\in\mathbb{C}:\ |z|=1 \}$.
By $\mathbb{Z}_{k,l}$ we mean all integers $r$, which satisfy the following inequality:
$k\leq r\leq l$.
By $\mathbb{P}$ we mean a set of all complex polynomials.
By $\mathbb{P}_n$ we denote a set of all complex polynomials, which have degrees less than or equal to $n$, $n\in\mathbb{Z}_+$.
By $\mathfrak{B}(M)$ we denote the set of all Borel subsets of $M\subseteq\mathbb{C}$.
For a measure $\mu$ on $\mathfrak{B}(M)$ we denote by $L^2_\mu = L^2_\mu(M)$ the usual space of all
(classes of equivalence of) Borel measurable complex-valued functions $f$ on $M$, such that
$\int_M |f|^2 d\mu < +\infty$.

If H is a Hilbert space then $(\cdot,\cdot)_H$ and $\| \cdot \|_H$ mean
the scalar product and the norm in $H$, respectively. 
Indices may be omitted in obvious cases.
All Hilbert spaces appearing in this paper are assumed to be separable.
For a linear operator $A$ in $H$, we denote by $D(A)$
its  domain, by $R(A)$ its range, and $A^*$ means the adjoint operator
if it exists. If $A$ is invertible then $A^{-1}$ means its
inverse. $\overline{A}$ means the closure of the operator, if the
operator is closable. If $A$ is bounded then $\| A \|$ denotes its
norm.
For a set $M\subseteq H$
we denote by $\overline{M}$ the closure of $M$ in the norm of $H$.
By $\mathop{\rm Lin}\nolimits M$ we mean
the set of all linear combinations of elements from $M$,
and $\mathop{\rm span}\nolimits M:= \overline{ \mathop{\rm Lin}\nolimits M }$.
By $E_H$ we denote the identity operator in $H$, i.e. $E_H x = x$,
$x\in H$. In obvious cases we may omit the index $H$. If $H_1$ is a subspace of $H$, then $P_{H_1} =
P_{H_1}^{H}$ denotes the orthogonal projection of $H$ onto $H_1$.

\section{Truncated moment problems on $\mathbb{C}$ and the similarity.}

At first we shall characterize those linear operators on a finite-dimensional Hilbert space $H$, which belong to the class $C_+(H)$.
We denote by
$\Gamma(y_0,y_1,...,y_n)$, the determinant of the Gram matrix of vectors $y_0,...,y_n\in H$, $n\in\mathbb{Z}_+$:
$$ \Gamma(y_0,y_1,...,y_n) = \det\left(
(y_k,y_l)_H
\right)_{k,l=0}^n. $$
The following theorem is an analog of Theorem~1 in~\cite{cit_3500_Z}.

\begin{theorem}
\label{t2_1}
Let $A$ be a linear operator in a $d$-dimensional Hilbert space $H$, $d>1$. The operator $A$ belongs to $C_+(H)$ if and only if
the following conditions hold:

\begin{itemize}
\item[(i)] $A$ is a cyclic complex symmetric operator in $H$;

\item[(ii)] there exists a cyclic vector $x_0$ of $A$ such that
\begin{equation}
\label{f2_5}
\Gamma(x_0,x_1,...,x_n,x_n^*) = 0,\qquad n\in\mathbb{Z}_{1,d-1},
\end{equation}
where
$$ x_k = A^k x_0,\quad x_k^* = (A^*)^k x_0, $$
and $Jx_0 = x_0$, for a conjugation $J$ in $H$, such that $JAJ=A^*$.

\end{itemize}

\end{theorem}
\textbf{Proof.} The proof goes along the same lines as the proof of Theorem~1 in~\cite{cit_3500_Z}, with some necessary modifications.
For convenience of the reader, we shall briefly present the arguments.

\noindent
\textit{Necessity.}
Let $\{ u_k \}_{k=0}^{d-1}$ be an orthonormal basis in $H$ such that 
$\mathcal{M} = ( (A u_l,u_k) )_{k,l=0}^{d-1} \in \mathfrak{M}_{d;3}^+$.
Observe that
$$ A u_k = m_{k-1,k} u_{k-1} + m_{k,k} u_k + m_{k+1,k} u_{k+1},\qquad k\in\mathbb{Z}_{0,d-1}, $$
where $m_{-1,0} = m_{d,d-1} := 0$.
Using this relation and the induction argument we get
\begin{equation}
\label{f2_10}
u_r \in \mathop{\rm Lin}\nolimits \left\{ u_0, A u_0,..., A^r u_0 \right\},\qquad r\in\mathbb{Z}_{0,d-1}.
\end{equation}
Therefore $H = \mathop{\rm Lin}\nolimits \{ A^k u_0 \}_{k=0}^{d-1}$, and $u_0$ is a cyclic vector of $A$.
Consider the following conjugation in $H$:
$$ J \sum_{k=0}^{d-1} \xi_k u_k = \sum_{k=0}^{d-1} \overline{\xi_k} u_k,\qquad \xi_k\in\mathbb{C}. $$
Since $\mathcal{M}$ is complex symmetric, then $A$ is a $J$-symmetric operator.
Denote $H_r = \mathop{\rm Lin}\nolimits \{ A^k u_0 \}_{k=0}^r$, $r\in\mathbb{Z}_{0,d-1}$.
By~(\ref{f2_10}) we see that $u_0,...,u_r\in H_r$, therefore they form an orthonormal basis in $H_r$.
Since $J u_k = u_k$, we get $J H_r\subseteq H_r$. Then
$$ (A^*)^r u_0 = (JAJ)^r u_0 = J A^r Ju_0 = J A^r u_0 \in H_r. $$
Vectors $u_0, A u_0, ..., A^r u_0, (A^*)^r u_0$, are linearly dependent, and we obtain relation~(\ref{f2_5}).

\noindent
\textit{Sufficiency.}
For a given cyclic vector $x_0$ we denote $H_r = \mathop{\rm Lin}\nolimits \{ A^k x_0 \}_{k=0}^r$, $r\in\mathbb{Z}_{0,d-1}$.
Notice that 
\begin{equation}
\label{f2_15}
A^{r+1} x_0 \notin H_r,\qquad r\in\mathbb{Z}_{0,d-2}.
\end{equation}
In fact, suppose to the contrary that $A^{r+j} x_0 \in H_r$, $1\leq j\leq k$, for some $r\in\mathbb{Z}_{0,d-2}$, $k\in\mathbb{N}$.
Then
$$ A^{r+k+1} x_0 = A A^{r+k} x_0 = A \sum_{t=0}^r \alpha_{r,k;t} A^t x_0 =
\sum_{t=0}^r \alpha_{r,k;t} A^{t+1} x_0\in H_r,\ (\alpha_{r,k;t}\in\mathbb{C}). $$
Repeating this trick we get $A^{r+1} x_0,A^{r+2} x_0, ....$,  all belong to $H_r$. Therefore $H=H_r$, a contradiction.
Applying the Gram-Schmidt orthogonalization process to 
$$ x_0, Ax_0,..., A^{d-1} x_0, $$
we get an orthonormal basis $\{ g_j \}_{j=0}^{d-1}$ in $H$. It also follows by the construction that
$\{ g_j \}_{j=0}^{r}$ is an orthonormal basis in $H_r$.
By~(\ref{f2_5}) we may write:
$$ J A^n x_0 = J A^n J x_0 = (A^*)^n x_0 \in H_n,\qquad n\in\mathbb{Z}_{1,d-1}. $$
Therefore $J H_r \subseteq H_r$. Let
$$ J g_r = \sum_{j=0}^r \beta_{r,j} g_j,\qquad \beta_{r,j}\in\mathbb{C},\quad r\in\mathbb{Z}_{0,d-1}. $$
Since $\beta_{r,k} = (Jg_r,g_k) = \overline{ (g_r,Jg_k) } = 0$, for $k\in\mathbb{Z}_{0,r-1}$, then
$J g_r = \beta_{r,r} g_r$.
Since $\| J g_r \| = \| g_r \| = 1$, then $\beta_{r,r} = e^{i\varphi_r}$, $\varphi_r\in[0,2\pi)$.
We set $u_r := e^{\frac{1}{2} \varphi_r i} g_r$, $r\in\mathbb{Z}_{0,d-1}$. Then $J u_r = u_r$.

Let us check that the matrix $\mathcal{M} = (m_{k,l})_{k,l=0}^{d-1} = ( (A u_l,u_k) )_{k,l=0}^{d-1}$, belongs to
$\mathfrak{M}_{d;3}^+$.
The complex symmetry of $\mathcal{M}$ follows from the complex symmetry of $A$. If $l>k+1$, then
$$ m_{k,l} = (A u_l, u_k) = (u_l, JA u_k) = 0, $$
since $J A u_k\in H_{k+1}$.
Notice that $A u_r \in H_{r+1}$, $r\in\mathbb{Z}_{0,d-2}$, since $u_r\in H_r$.
Observe that $A u_r \notin H_r$, $r\in\mathbb{Z}_{0,d-2}$. In the opposite case we would get $A H_r\subseteq H_r$. Then
$A^k x_0 \subseteq H_r$, $k\in\mathbb{Z}_+$, and $H=H_r$, a contradiction.
Therefore for $r\in\mathbb{Z}_{0,d-2}$, we may write:
$$ A u_r = \sum_{j=0}^{r+1} \gamma_{r,j} u_j,\qquad \gamma_{r,j}\in\mathbb{C},\ \gamma_{r,r+1}\not=0, $$
and $m_{r+1,r} = \gamma_{r,r+1} \not=0$.
$\Box$

We shall now turn to the study of the moment problem~(\ref{f1_50}).
Let $\{ s_k \}_{k\in\mathbb{Z}_{0,\rho}}$ be a prescribed set of complex numbers, such that $s_0 > 0$, ($\rho\in\mathbb{N}$).
If $\rho = 1$, it is seen that a $1$-atomic measure with an atom at $\frac{s_1}{s_0}$,
having a mass $s_0$, is a solution of the moment problem~(\ref{f1_50}).
Thus, we shall assume that $\rho \geq 2$.

\begin{lemma}
\label{l2_1}
Let the moment problem~(\ref{f1_50}) be given, with $\rho\geq 2$, and the following given moments:
$$ s_0 = 1,\ s_\rho = c,\quad s_k = 0,\quad k\in\mathbb{Z}_{1,\rho-1}, $$
where $c$ is an arbitrary complex number.
Then the moment problem~(\ref{f1_50}) has a finitely atomic solution with atoms, situated on the circle
$T_r := \{ z\in\mathbb{C}:\ |z| = r \}$, where $r$ is an arbitrary positive number, greater than 
$\sqrt[\rho]{|c|}$.
\end{lemma}
\textbf{Proof.}
Denote $\widetilde c := \frac{ c }{ r^\rho }$. Notice that $|\widetilde c| < 1$.
Observe that the determinant of the following Toeplitz matrix of size $(\rho+1)\times(\rho+1)$:
$$ T_\rho =
\left(
\begin{array}{ccccc}
1 & 0 & \ldots & 0 & \widetilde{c} \\
0 & 1 & \ldots & 0 & 0 \\
\vdots & \vdots & \ddots & \vdots & \vdots \\
0 & 0 & \ldots & 1 & 0 \\
\overline{\widetilde c} & 0 & \ldots & 0 & 1 \end{array}
\right), $$
is equal to $1 - |\widetilde{c}|^2 > 0$.
Therefore the truncated trigonometric moment problem with moments (see, e.g.,~\cite{cit_2000__Krein_Nudelman}):
$$ s_0' = 1,\ s_\rho' = \widetilde{c},\quad s_k' = 0,\quad k\in\mathbb{Z}_{1,\rho-1}, $$
has a finitely atomic solution:
$$ \sum_j z_j^k m_j = s_k' = \frac{s_k}{r^k},\qquad k\in\mathbb{Z}_{0,\rho}, $$
where $z_j\in\mathbb{T}$, $m_j > 0$.
Then
$$ \sum_j u_j^k m_j = s_k,\qquad k\in\mathbb{Z}_{0,\rho}, $$
where $u_j := r z_j\in T_r$.
$\Box$

Thus, we can construct finitely atomic solutions of the moment problem~(\ref{f1_50}), with $\rho = 1$, $s_0 > 0$, and with
$\rho\geq 2$, when $s_0 > 0$, $s_\rho\in\mathbb{C}$, $s_k = 0$, $k\in\mathbb{Z}_{1,\rho-1}$.
In fact, one can divide the moments by $s_0$ and apply Lemma~\ref{l2_1}.
We are ready to present an algorithm for the moment problem~(\ref{f1_50}).

\noindent
\textbf{Algorithm 1.}

\noindent
\textbf{Input data.}
$\rho\in\mathbb{Z}:\ \rho\geq 2$.
Moments $s_0 > 0$, $s_k\in\mathbb{C}$, $k\in\mathbb{Z}_{1,\rho}$.

\noindent
\textbf{Step 1.} Construct a solution $\mu_1$ to the moment problem~(\ref{f1_50}),
with the following moments:
$$ s_0(\mu_1) = \frac{ s_0 }{ \rho },\quad s_1(\mu_1) = s_1, $$
see considerations before Lemma~\ref{l2_1}.

\noindent
\textbf{Step n, $n=2,...,\rho$.} Construct a solution $\mu_n$ to the moment problem~(\ref{f1_50}), with the following moments:
$$ s_0(\mu_n) = \frac{ s_0 }{ \rho },\quad s_n(\mu_n) = s_n - \sum_{l=1}^{n-1} s_n(\mu_l),\quad  s_j(\mu_n) = 0,\ j\in\mathbb{Z}_{1,n-1}, $$
see the proof of Lemma~\ref{l2_1}.

\noindent
\textbf{Step $\rho+1$.} Set
$$ \mu = \sum_{l=1}^\rho \mu_l. $$

\noindent
\textbf{Output.}
A finitely atomic solution $\mu$ of the moment problem~(\ref{f1_50}).

Let us illustrate the above algorithm by a numerical example.

\begin{example}
\label{e2_1}
Consider the moment problem~(\ref{f1_50}), with the following moments:
$$ s_0 = 1,\quad s_1 = 1+i,\quad s_2 = 3i. $$

\noindent
\textbf{Step 1}. We need a solution $\mu_1$ of the moment problem~(\ref{f1_50}), with the following moments:
$$ s_0(\mu_1) = \frac{1}{2},\quad s_1(\mu_1) = 1+i. $$
One can choose $\mu_1$ to be
a $1$-atomic measure with an atom at $2+2i$,
with a mass $\frac{1}{2}$.

\noindent
\textbf{Step 2}. We now need a solution $\mu_2$ of the moment problem~(\ref{f1_50}), with the following moments:
\begin{equation}
\label{f2_17}
s_0(\mu_2) = \frac{1}{2},\quad s_1(\mu_2) = 0,\quad s_2(\mu_2) = s_2 - s_2(\mu_1) = -i. 
\end{equation}
Let us construct a solution to the moment problem~(\ref{f1_50}), with the normalized moments:
\begin{equation}
\label{f2_18}
\widehat s_0 = 1,\quad \widehat s_1 = 0,\quad \widehat s_2 = -2i. 
\end{equation}
According to Lemma~\ref{l2_1}, we choose $r = 2$.
Following the construction in the proof of Lemma~1, we set $\widetilde c = -\frac{i}{2}$.
Then the corresponding truncated trigonometric moment problem with moments
\begin{equation}
\label{f2_20}
s_0' = 1,\quad s_1' = 0,\quad s_2' = -\frac{i}{2}, 
\end{equation}
is solvable.
In order to construct a solution we may use known descriptions. For example, one can use an operator description 
from~\cite{cit_3550_Z_TTMP}.
We obtain that the truncated trigonometric moment problem~(\ref{f2_20}) has a $3$-atomic solution with atoms
at 
$$ z_0 = \frac{ 1 }{ \sqrt{2} } (1+i),\quad z_{1,2} = \frac{1}{ 4\sqrt{2} } \left(
-1\mp\sqrt{15} + i ( -1\pm\sqrt{15} )
\right), $$
and masses
$m_0=\frac{1}{5}$, $m_1=m_2=\frac{2}{5}$.
Then the moment problem~(\ref{f1_50}) with the normalized moments~(\ref{f2_18}) has 
a $3$-atomic solution with atoms at $2z_j$, and masses $m_j$.
Consequently, we can choose $\mu_2$ to be a $3$-atomic solution with atoms at $2z_j$, and masses $\frac{1}{2} m_j$, $j=0,1,2$.

\noindent
\textbf{Step 3}. We set $\mu = \mu_1 + \mu_2$.
Thus, $\mu$ is a $4$-atomic solution with 
an atom at $2+2i$, with a mass $\frac{1}{2}$, and
atoms at $2z_j$, and masses $\frac{1}{2} m_j$, $j=0,1,2$.

\end{example}

\begin{remark}
\label{r2_1}
It would be of interest to adapt Algorithm~1 for the full moment problem consisting of finding a measure $\mu$ on $\mathfrak{B}(\mathbb{C})$
such that relation~(\ref{f1_50}) holds for $k\in\mathbb{Z}_+$, with a set
$\{ s_k \}_{k\in\mathbb{Z}_+}$ of given complex numbers.
One can consider similar steps but with $s_0(\mu_n) = \frac{1}{2^{n+1}} s_0$, $n=1,2,...$.
However, it is not clear whether $\sum_{n=0}^\infty \mu_n$ has absolutely convergent moments.

\end{remark}

We can now state the main result.
\begin{theorem}
\label{t2_2}
Let $A$ be a linear operator in a $d$-dimensional Hilbert space $H$, $d>1$. If the operator $A$ belongs to $C_+(H)$, 
then it is similar to a rank-one perturbation of a restriction of a normal operator.
Namely, there exists a finitely atomic measure $\mu$ on $\mathfrak{B}(\mathbb{C})$, and an invertible linear operator
$T$, which maps $H$ into $L^2_\mu$, such that:
\begin{equation}
\label{f2_25}
T A T^{-1} = Z_0 + a(z) (\cdot, b(z))_{L^2_\mu},
\end{equation}
where $Z_0$ is the operator of the multiplication by an independent variable in $L^2_\mu$, restricted
to a linear set of all complex polynomials with degrees less than or equal to $d-1$;
$a(z)$ and $b(z)$ are some complex polynomials of $z$ and $\overline{z}$, respectively.

\end{theorem}
\textbf{Proof.} 
Since the given operator $A$ belongs to $C_+(H)$, then there exists an
orthonormal basis $\{ u_k \}_{k=0}^{d-1}$ in $H$ such that the matrix $\mathcal{M}$ from~(\ref{f1_35})
belongs to $\mathfrak{M}_{d;3}^+$.
Thus $\mathcal{M} = (m_{k.l})_{k.l=0}^{d-1}$ is a three-diagonal complex symmetric matrix, with non-zero entries on the
first sub-diagonal. Let us extend $\mathcal{M}$ to a semi-infinite matrix $J = (m_{k.l})_{k.l=0}^\infty$, setting
$$ m_{k+1,k} = m_{k,k+1} = 1,\qquad k=d-1,d,..., $$
ant setting $m_{k,l}=0$, for the rest of new entries.
Then $J$ is a three-diagonal semi-infinite complex symmetric matrix.
Denote 
$$ a_k := m_{k,k+1},\qquad b_k := m_{k,k},\quad k\in\mathbb{Z}_+. $$
According to the procedure in~\cite{cit_2500_Z_Serdica}, we define a sequence of polynomials
$\{ p_n(\lambda) \}_{n=0}^\infty$, $p_0(\lambda)=1$, by the following recurrence relations:
$$ b_0 p_0(\lambda) + a_0 p_1(\lambda) = \lambda p_0(\lambda), $$
\begin{equation}
\label{f2_30}
a_{n-1} p_{n-1}(\lambda) + b_n p_n(\lambda) + a_n p_{n+1}(\lambda) = \lambda p_n(\lambda),\qquad n\in\mathbb{N}.
\end{equation}
Then we define \textit{a linear with respect to both arguments} functional $\sigma(u,v)$, $u,v\in\mathbb{P}$, (the spectral function)
by the following relation:
\begin{equation}
\label{f2_34}
\sigma(p_n(\lambda),p_m(\lambda)) = \delta_{m,n},\qquad m,n\in\mathbb{Z}_+, 
\end{equation}
and extending by the linearity.
This functional obeys the following property (see formula~(10) in~\cite{cit_2500_Z_Serdica}):
\begin{equation}
\label{f2_35}
\sigma(\lambda u(\lambda),v(\lambda)) = \sigma(u(\lambda),\lambda v(\lambda)),\qquad u,v\in\mathbb{P}. 
\end{equation}
Define a linear functional $S(u)$, $u\in\mathbb{P}$, in the following way:
$$ S(u(\lambda)) = \sigma(u(\lambda),1),\qquad u\in\mathbb{P}. $$
By~(\ref{f2_35}) it is seen that
\begin{equation}
\label{f2_36}
S(u(\lambda)v(\lambda)) = \sigma(u(\lambda),v(\lambda)),\qquad u,v\in\mathbb{P}. 
\end{equation}
Denote
\begin{equation}
\label{f2_38}
s_k = S(\lambda^k),\qquad k\in\mathbb{Z}_{0,\rho},\qquad  \rho\in\mathbb{Z}: \rho>2d. 
\end{equation}
Notice that $s_0 = \sigma(p_0 p_0) = 1$.
By~Algorithm~1 we can construct a finitely atomic solution $\mu$ to the moment problem~(\ref{f1_50}) with moments~(\ref{f2_38}):
$$ \int_{\mathbb{C}} z^k d\mu = s_k,\qquad k\in\mathbb{Z}_{0,\rho}. $$
By~(\ref{f2_36}),(\ref{f2_34}) we obtain that
\begin{equation}
\label{f2_42}
\int_{\mathbb{C}} p_n(\lambda) p_m(\lambda) d\mu = \delta_{m,n},\qquad m,n\in\mathbb{Z}_{0,d}. 
\end{equation}
Let us check that
\begin{equation}
\label{f2_44}
u_k = p_k(A) u_0,\qquad k\in\mathbb{Z}_{0,d-1}. 
\end{equation}
We shall proceed by the induction argument.
If $k=0$, then relation~(\ref{f2_44}) holds.
Suppose that relation~(\ref{f2_44}) holds for $k\in\mathbb{Z}_{0,n}$, with some $n\in\mathbb{Z}_{0,d-2}$.
We need to check it for $k=n+1$. By~(\ref{f2_30}) we may write:
$$ p_{n+1}(A) = \frac{1}{a_n} (A p_n(A) - a_{n-1} p_{n-1}(A) - b_n p_n(A) ), $$
where $p_{-1}:=0$, $a_{-1}:= 0$.
Applying both sides to $u_0$ we get:
$$ p_{n+1}(A) u_0 = \frac{1}{a_n} (A p_n(A) u_0 - a_{n-1} p_{n-1}(A) u_0 - b_n p_n(A) u_0) = $$
$$ = \frac{1}{a_n} (A u_n - a_{n-1} u_{n-1} - b_n u_n) = u_n, $$
since $A u_n = m_{n-1,n} u_{n-1} + m_{n,n} u_n + m_{n+1,n} u_{n+1} = a_{n-1} u_{n-1} + b_n u_n + a_n u_{n+1}$.
Thus, relation~(\ref{f2_44}) holds true. Here the complex symmetry of $\mathcal{M}$ played an essential role,
since the polynomial $p_n$ was defined by the coefficients in the $n$-th row of $J$ (as usual, a vector of polynomials forms a generalized
eigenvalue of $J$), while $A u_n$ used coefficients
from the $n$-th column of $\mathcal{M}$.
Define the following linear operator $T$, which maps $H$ into $L^2_\mu$:
\begin{equation}
\label{f2_45}
T \sum_{k=0}^{d-1} \xi_k u_k = \sum_{k=0}^{d-1} \xi_k p_k(z),\qquad \xi_k\in\mathbb{C}. 
\end{equation}
Observe that in the $k$-th step of Algorithm~1 one constructs a measure $\mu_k$, with atoms on a circle $T_r$. This circle
can be chosen to have an arbitrarily big radius $r$. Thus, we can assume that all these circles on different steps
do not intersect. Since there is at least one atom on each circle, then the number of atoms is bigger than or equal to $\rho$.
Consequently, in the case of the moments~(\ref{f2_38}), we can assume that we have a solution $\mu$ with more then $2d$ atoms.
The latter fact implies that the operator $T$ is invertible. In fact, suppose that
$$ T u = 0,\qquad u = \sum_{k=0}^{d-1} \xi_k u_k\in H,\ \xi_k\in\mathbb{C}. $$
Then
$$ 0 = \left\| \sum_{k=0}^{d-1} \xi_k p_k(z) \right\|_{L^2_\mu}^2 =
\int_\mathbb{C} \left|
\sum_{k=0}^{d-1} \xi_k p_k(z)
\right|^2 d\mu. $$
Therefore the polynomial $\sum_{k=0}^{d-1} \xi_k p_k(z)$, of degree less that or equal to $d-1$, should vanish at $2d$ points.
Thus, we get $\sum_{k=0}^{d-1} \xi_k p_k(z) = 0$. Since $\deg p_k = k$, it follows that all $\xi_k$ are zero, and $u=0$.

Denote by $L_0$ a linear set in $L^2_\mu$ of all (classes of the equivalence containing) complex polynomials with degrees less that $d$.
Observe that polynomials $\{ p_k(z) \}_{k=0}^{d-1}$ are linearly independent and they form a linear basis in $L_0$.
Let $u(z)$ be an arbitrary element of $L_0$:
$$ u(z) = \sum_{j=0}^{d-2} \alpha_j p_j(z) + c p_{d-1}(z),\qquad  \alpha_j,c\in\mathbb{C}. $$
By~(\ref{f2_30}),(\ref{f2_42}) we may write
$$ TAT^{-1} u(z) = TA 
\left(
\sum_{j=0}^{d-2} \alpha_j u_j + c u_{d-1}
\right) = $$
$$ = T\left( \sum_{j=0}^{d-2} \alpha_j ( a_{j-1} u_{j-1} + b_j u_j + a_j u_{j+1} )
+ c ( a_{d-2} u_{d-2} + b_{d-1} u_{d-1} ) 
\right) = $$
$$ = \sum_{j=0}^{d-2} \alpha_j ( a_{j-1} p_{j-1}(z) + b_j p_j(z) + a_j p_{j+1}(z) )
+ c ( a_{d-2} p_{d-2}(z) + b_{d-1} p_{d-1}(z) )  = $$
$$ = z \left( \sum_{j=0}^{d-2} \alpha_j p_j(z) + c p_{d-1}(z) \right) - a_{d-1} c p_d(z) = $$
$$ = Z_0 u(z) - a_{d-1} p_d(z) (u(z),\overline{ p_{d-1}(z) })_{L^2_\mu}, $$
and relation~(\ref{f2_25}) is proved.
$\Box$

The following open problems seems to be interesting for further investigations:

\noindent
\textbf{Problem 1.} \textit{Describe those complex symmetric operators which are similar to normal operators}.

\noindent
\textbf{Problem 2.} \textit{Describe those complex symmetric operators which are similar to self-adjoint / unitary operators}.

\noindent
\textbf{Problem 3.} \textit{Describe those complex symmetric operators which are similar to finite-rank perturbations of normal operators}.

\noindent
\textbf{Problem 4.} \textit{Solve the full moment problem from Remark~\ref{r2_1}}.

We hope that all these questions will be solved in the near future.

}

\noindent
Address:

V. N. Karazin Kharkiv National University \newline\indent
School of Mathematics and Computer Sciences \newline\indent
Department of Higher Mathematics and Informatics \newline\indent
Svobody Square 4, 61022, Kharkiv, Ukraine

Sergey.M.Zagorodnyuk@gmail.com; Sergey.M.Zagorodnyuk@univer.kharkov.ua

\end{document}